\documentclass[oneside,notitlepage,12pt]{article}

\usepackage{amssymb}
\usepackage[leqno]{amsmath}
\usepackage{amsfonts}
\usepackage{amsopn}
\usepackage{amstext}
\usepackage{amsthm}

\usepackage{enumitem}

\usepackage{tikz}
\usetikzlibrary{cd}

\usepackage{verbatim}
\usepackage[colorlinks]{hyperref}
\usepackage{makeidx}

\newcommand{\define}[2]{{\em #1}\index{#2}}

\usepackage{calrsfs}
\usepackage{fourier-orns}
\usepackage{clock} 
\usepackage[alpine, weather]{ifsym}

\providecommand{\cal}{\mathcal}
\renewcommand{\Bbb}{\mathbb}

\newenvironment{pf}{\begin{proof}}{\end{proof}}

\newcommand{\Ef}{{\cal{F}}}
\newcommand{\Gee}{{\cal{G}}}

\newcommand{\El}{{\cal{L}}}
\newcommand{\Pee}{{\cal{P}}}

\newcommand{\Vee}{{\cal{V}}}
\newcommand{\Wu}{{\cal{W}}}
\newcommand{\Zee}{{\Bbb{Z}}}

\newcommand{\Qyu}{{\Bbb{Q}}}

\newcommand{\sig}{\sigma}

\renewcommand{\phi}{\varphi}
\renewcommand{\rho}{\varrho}

\newcommand{\rest}{\restriction}

\newcommand{\ntr}{{n\in\omega}}

\newcommand{\loe}{\leq}

\newcommand{\subs}{\subseteq}
\newcommand{\sups}{\supseteq}
\newcommand{\nnempty}{\ne\emptyset}
\newcommand{\argum}{\:\cdot\:}

\renewcommand{\iff}{\Longleftrightarrow}

\newcommand{\id}[1]{{\operatorname{i\!d}_{#1}}} 

\newcommand{\symd}{\div} 

\newtheorem{thm}{Theorem} 
\newtheorem{tw}{Theorem}[section]
\newtheorem{wn}[tw]{Corollary}
\newtheorem{lm}[tw]{Lemma}
\newtheorem{prop}[tw]{Proposition}

\theoremstyle{definition}

\theoremstyle{remark}

\newcommand{\setof}[2]{\{#1\colon #2\}}

\newcommand{\seq}[1]{\langle #1 \rangle}

\newcommand{\sett}[2]{\{#1\}_{#2}}
\newcommand{\sn}[1]{\{#1\}} 
\newcommand{\dn}[2]{\{#1,#2\}} 
\newcommand{\pair}[2]{\langle #1, #2 \rangle} 
\newcommand{\triple}[3]{\langle #1, #2, #3 \rangle} 
\newcommand{\map}[3]{#1\colon #2 \to #3} 
\newcommand{\img}[2]{#1[#2]} 

\newcommand{\fra}{Fra\"iss\'e}

\newcommand{\U}{\mathbb U}

\providecommand{\nat}{\omega}

\newcommand{\ciag}[1]{{\sett{{#1}_n}{\ntr}}}

\newcommand{\aut}{\operatorname{Aut}}

\newcommand{\iso}{\approx}

\newcommand{\bX}{{\mathbb{X}}}

\newcommand{\cmp}{\circ} 

\newcommand{\bW}{\mathbb W}

\newcommand{\loc}{\sigma} 

\newcommand{\procons}{\mathcal{C\!P}}
\newcommand{\kola}{{\mathcal{Z}}}

\newcommand{\dupg}[2]{{\mathbb D}^{#1}#2}
\newcommand{\dup}[1]{\dupg{}#1}
\newcommand{\cR}{\mathcal R}
\newcommand{\suk}[1]{\operatorname{s}(#1)}
\newcommand{\misa}{\pm}

\title{Homogeneous structures with non-universal automorphism groups}
\author{
{\sc Wies{\l}aw Kubi\'s}\footnote{Supported by GA\v CR grant No. 17-27844S (Czech Science Foundation).}\\
{\small Mathematical Institute, Czech Academy of Sciences, Czech Republic}\\
{\small Cardinal Stefan Wyszy\'nski University, Warsaw, Poland}
\and
{\sc Saharon Shelah}\footnote{Partially supported by European Research Council grant 338821. Paper 1153 on Shelah's list.}\\
{\small Hebrew University of Jerusalem, Israel}\\
{\small Rutgers University, USA}
}

\date{\clocktime\today}

\begin{document}

\maketitle

\begin{abstract}
We present three examples of countable homogeneous structures (also called \emph{\fra\ limits}) whose automorphism groups are not universal, namely, fail to contain isomorphic copies of all automorphism groups of their substructures.

Our first example is a particular case of a rather general construction on \fra\ classes, which we call \emph{diversification}, leading to automorphism groups containing copies of all finite groups.
Our second example is a special case of another general construction on \fra\ classes, the \emph{mixed sums}, leading to a \fra\ class with all finite symmetric groups appearing as automorphism groups and at the same time with a torsion-free automorphism group of its \fra\ limit.
Our last example is a \fra\ class of finite models with arbitrarily large finite abelian automorphism groups, such that the automorphism group of its \fra\ limit is again torsion-free.

\


\

\end{abstract}

\newpage

\tableofcontents

\section{Introduction}

This note concerns mathematical structures with a high level of symmetry. Symmetries are automorphisms, namely, bijections fully preserving the structure. An object is often called \emph{homogeneous} if ``small" isomorphisms between its sub-objects extend to automorphisms. The meaning of ``small" depends on the context. We are interested in countable models of a first order language, where ``small" means ``finite". Specifically, our objects of study are \emph{\fra\ limits}, well known in model theory since the work of \fra~\cite{Fraisse}. In this setting, being homogeneous means, in particular, that every automorphism between finite submodels extends to a full automorphism. In most of the well known concrete cases such an extension can be made uniform, in the sense that it preserves composition and yields an embedding between the automorphism groups. In particular, the automorphism group of such a \fra\ limit contains isomorphic copies of all automorphism groups of its finite submodels. A homogeneous structure admitting extension operators preserving compositions can be called \emph{uniformly homogeneous}.
As it happens, not all \fra\ limits are uniformly homogeneous. The purpose of this note is to provide suitable counterexamples.

\paragraph{}
Let $\Ef$ be a \fra\ class and let $\loc \Ef$ denote the class of all models $M = \bigcup_{\ntr} E_n$, where $\ciag{E}$ is a chain in $\Ef$.
Let $U$ be the \fra\ limit of $\Ef$.
We say that $\aut(U)$ is \emph{universal} if for every $M \in \loc \Ef$ the group $\aut(M)$ is isomorphic to a subgroup of $\aut(U)$.
Most of the natural and typical \fra\ classes have this property. In fact, it is guaranteed by the existence of so-called Kat\v etov functors, see~\cite{KubMas}.
In fact, in the presence of a Kat\v etov functor, the group $\aut(U)$, where $U$ is the \fra\ limit, is universal also in the topological sense (recall that $\aut(M)$ carries a natural Polish topology for every countable model $M$).
The question whether $\aut(U)$ is always universal was asked by Eric Jaligot~\cite{Jaligot} and perhaps also by some other people working in this topic. In particular, the question is repeated recently in \cite{Doucha} and \cite{Isabel}; the latter paper contains a remark (attributed to Piotr Kowalski) that the class of all finite fields of a fixed characteristic provides a counterexample. Nevertheless, the question remained open when asking for relational classes.
We answer it in the negative.
Namely, we prove:

\begin{thm}\label{ThmMejnOne}
	There exists a relational \fra\ class in a finite signature, such that the group of automorphisms of its \fra\ limit is not universal, whereas it contains isomorphic copies of all finite groups.
\end{thm}

\begin{thm}\label{ThmMejnTwo}
	There exists a relational \fra\ class in a finite signature, such that the automorphism group of its \fra\ limit is torsion-free, while
	the class of automorphism groups of its finite substructures contains all finite symmetric groups.
\end{thm}

\begin{thm}\label{ThmMejnTwoAnAhalf}
	There exists a \fra\ class of finite models $\kola$ in a finite signature, such that the class $$\setof{\aut(M)}{M \in \kola}$$ consists of abelian groups, contains all possible finite products of finite cyclic groups, while the automorphism group of its \fra\ limit is torsion-free.
\end{thm}

In the last result the signature consists of two binary relation symbols and one unary function symbol, so it is not relational. The first result shows that the automorphism group of a relational \fra\ limit can possibly contain isomorphic copies of all automorphism groups of its age, while still being non-universal in the sense described above.
The second and the third results show the more extreme situations, where the automorphism group of a homogeneous structure does not contain any non-trivial finite groups.

\section{Preliminaries}

We shall use standard notation concerning model theory. For undefined notions we refer to~\cite{Hodges}.
Recall that a class $\Ef$ of finitely generated models is a \emph{\fra\ class} if it is hereditary, has the joint embedding property, countably many isomorphic types, and the amalgamation property.
The \emph{amalgamation property} says that for every two embeddings $\map f Z X$, $\map g Z Y$ there exist $W$ and embeddings $\map {f'} X W$, $\map {g'} Y W$ such that $f' f = g' g$.
If one can always have $\img {f'}X \cap \img {g'}Y = \img{f' f}{Z}$ then this is called the \emph{strong amalgamation property}. Replacing $f$ and $g$ by inclusions and assuming $Z = X \cap Y$, the strong amalgamation property means that one can amalgamate $X$ and $Y$ over $Z$ without identifying points in $X \cup Y$.

Recall that every \fra\ class $\Ef$ has its unique \emph{\fra\ limit} $\U \in \sig \Ef$ which is characterized by the \emph{extension property}, namely, given structures $A \subs B$ in $\Ef$, every embedding of $A$ into $\U$ extends to a embedding of $B$. One needs to assume  also that $\Ef$ is the \emph{age} of $\U$, namely, $\Ef$ equals the class of all finitely generated structures embeddable into $\U$.
Of course, $\U$ is always homogeneous and conversely, every countable homogeneous structure is the \fra\ limit of its age. See~\cite{Hodges} for more details.

In this note we are interested in \fra\ classes of finite models. A class is \emph{relational} if its signature consists of relation symbols only.
In particular, the empty set is a model of any relational signature, since there are no constants.

We denote by $S_\infty$ the symmetric group (i.e. the group of all permutations) of a countable infinite set. $S_n$ will denote the symmetric group of the set $n = \{0,1, \dots, n-1\}$. The set of all natural numbers (including zero) will be denoted, as usual, by $\nat$.

\subsection{Some very basic group theory}

Given a set $A$, its powerset $\Pee(A)$ has a natural abelian group operation, namely, the symmetric difference $\symd$ (recall that $A \symd B = (A \setminus B) \cup (B \setminus A)$).
Note that each element of $\Pee(A)$ has degree $2$ (except its unit $\emptyset$).
The $2$-element group $\pair{\Pee(1)}{ \symd}$ is usually denoted by $\Zee_2$.

\begin{prop}
	Assume $G$ is a group with a subgroup $K$ of index $<2^{\aleph_0}$. If no element of $K$ has order $2$ then $\pair{\Pee(\omega)}\symd$ does not embed into $G$.
\end{prop}

\begin{pf}
	Suppose $A \mapsto f_A$ is an embedding of $\Pee(\omega)$ into $G$, that is, $f_A f_B = f_{A \symd B}$ for every $A, B \subs \omega$ and $f_A \ne 1$ for every $A \nnempty$.
	Since $|\Pee(\omega)| = 2^{\aleph_0}$, there are $A \ne B$ such that $f_A$ and $f_B$ belong to the same right co-set of $K$. Now
	$$f_{A \symd B} = f_A f_B = f_A (f_B)^{-1} \in K$$
	and $f_{A \symd B} \ne 1$ has order $2$, a contradiction.	
\end{pf}

Note that $\pair{\Pee(\nat)} \symd$ embeds into $S_\infty$. Namely, given $A \subs \nat \setminus 1$, define $\map{h_A}{\Zee}{\Zee}$ by $h_A(x) = -x$ if $x \in A$ and $h_A(x) = x$ otherwise.
Then $A \mapsto h_A$ is an embedding of $\Pee(\omega \setminus 1)$ into $S_\infty$. Thus:

\begin{wn}\label{CorWNebwgibEIG}
	Assume $S_\infty$ embeds into a group $G$. Then every subgroup of $G$ of index $<2^{\aleph_0}$ contains an element of order $2$.
\end{wn}

\subsection{Bipartite graphs}

We fix the notation concerning bipartite graphs.
Namely, these are structures of the form $\bX = \seq{X, L^\bX, R^\bX, \sim^\bX}$, where $L^\bX$, $R^\bX$ are unary predicates, $\sim^\bX$ is a symmetric binary relation, $\dn{L^\bX}{R^\bX}$ is a partition of $X$, and $x \sim^\bX y$ holds only if either $x \in L^\bX$, $y \in R^\bX$ or $x \in R^\bX$, $y \in L^\bX$. By this way, $\pair{X}{\sim^\bX}$ is indeed a bipartite graph and $L^\bX$, $R^\bX$ specifies its bipartition.
Adding the unary predicates $L$, $R$ to the signature, we fix the two sides of the bipartite graph, the one specified by $L$ could be called the \emph{left-hand side} while the other one could be called the \emph{right-hand side}. Note that
embeddings preserve the sides.
The following fact belongs to the folklore of \fra\ theory.

\begin{prop}\label{PROPbhvbivbw}
	The class of all finite bipartite graphs, described as above, is a \fra\ class. Its \fra\ limit is the unique countable bipartite graph
	$$\U = \seq{U, L^\U, R^\U, \sim^\U}$$
	satisfying the following condition:
	\begin{enumerate}
		\item[$(\star)$] For every finite disjoint sets $A, B \subs U$, there are $\ell \in L^\U \setminus (A \cup B)$, $r \in R^\U \setminus (A \cup B)$ such that $\ell \sim^\U y$ holds for $y \in A \cap R^\U$, $\ell \not \sim^\U y$ holds for $y \in B \cap R^\U$, $x \sim^\U r$ holds for $x \in A \cap L^\U$, and $x \not\sim^\U r$ holds for $x \in B \cap L^\U$.
	\end{enumerate}
\end{prop}

The structure $\U$ described above is called the \emph{universal homogeneous bipartite graph}.

\section{Diversifications}

In this section we present a general construction on \fra\ classes, leading to automorphism groups containing copies of all finite groups.

Let $\cR$ be a fixed countable relational signature (i.e. it consists of at most countably many relation symbols, no function symbols and no constant symbols).
Let $\Ef$ be a class of finite $\cR$-models.
We define the \define{diversification}{diversification} of $\Ef$, denoted by $\dup \Ef$, to be the class of two-sorted models $X = P^X \cup C^X$ with $P^X \cap C^X = \emptyset$ such that for each $y \in C^X$ the set $P^X$ is endowed with an $\cR$-structure $\cR_y$.
Clearly, this can be formalized in a first order language.
Specifically, for each $n$-ary relation symbol $R$ let $\tilde R$ be an $(n+1)$-ary relation symbol and let
$$\tilde \cR = \setof{\tilde R}{R \in \cR} \cup \dn{P}{C}.$$
Now, $\dup \Ef$ consists of finite models of $\tilde \cR$ satisfying the obvious axioms: $\dn P C$ is a partition and for each $y \in C$ the set $P$ endowed with $\cR_y := \setof{R( - , y)}{R \in \tilde \cR}$ is in $\Ef$.
Note that if $\Ef$ is hereditary then so is $\dup \Ef$, if $\Ef$ has countably many isomorphism types then so does $\dup \Ef$.
We shall see in a moment that if $\Ef$ is a \fra\ class with strong amalgamations then so is $\dup \Ef$.
Let us note that a similar construction, called \emph{generic variation}, is due to Baudisch~\cite{Baudisch}. The difference is that in a generic variation there are no extra predicates $P$ and $C$, hence $R(-,y)$ is defined for every $y$ and for every relation $R$ in the signature.

In order to avoid repetitions, we shall now introduce a more general version of diversifications, involving a group action.
Namely, fix a finite group $G$ and define $\dupg G \Ef$ to be the class of all models from $\dup \Ef$ with a distinguished free $G$-action, denoted by $\pair x g \mapsto x^g$.
More precisely, $\pair X a \in \dupg G \Ef$ if 
$X \in \dup \Ef$ and $\map a {X \times G} X$ is a free group action on $X$, that is, $\img a {C^X \times G} = C^X$, $\img a {P^X \times G} = P^X$ and, denoting $a(x,g) = x^g$, for every relation $R \in \cR$ the following implication holds.
\begin{equation}
	R(x_1,\dots, x_n, y) \implies R(x_1^g, \dots, x_n^g, y^g).
	\tag{$\pitchfork$}\label{EqPiczfork}
\end{equation}
In fact, the implication above is an equivalence.
Recall that a group action is \emph{free} if $x^g = x^h \implies g = h$ for every $x \in X$, $g,h \in G$.
Equivalently, if $x^g = x$ for some $x \in X$ then $g = 1$.
We are now ready to prove the crucial lemma.


\begin{lm}\label{LMKrucjalni}
	Assume $G$ is a finite group and $\Ef$ is a relational \fra\ class with strong amalgamations.
	Then $\dupg G \Ef$ is a \fra\ class.
\end{lm}

\begin{pf}
	First of all, note that either $\Ef = \sn \emptyset$ or else $\Ef$ contains arbitrarily large finite models, thanks to strong amalgamations.
	In the latter case, every model in $\Ef$ can be extended to an arbitrarily large finite model in $\Ef$, thanks to the joint embedding property.
	It is clear that $\dupg G \Ef$ is hereditary and has countably many isomorphism types.
	As the empty model is in $\Ef$, the joint embedding property of $\dupg G \Ef$ follows from the amalgamation property, which we prove below.

	Fix $Z \in \dupg G \Ef$.
	There are two types of ``simple" extensions of $Z$: adding a single $G$-orbit to $P$ and adding a single $G$-orbit to $C$. Clearly, all other embeddings are finite compositions of ``simple" ones.
	So let us fix ``simple" extensions $Z \subs X$, $Z \subs Y$. Without loss of generality, we may assume that $X \cap Y = Z$.
	Let $W = X \cup Y$. Then $W$ already has a uniquely determined free $G$-action. It remains to define a $\dup \Ef$-structure on $W$, so that the $G$-action will consist of automorphisms of $W$, i.e., that condition (\ref{EqPiczfork}) holds true.
	
	Define $P^W = P^X \cup P^Y$, $C^W = C^X \cup C^Y$. Clearly, $\dn{P^W}{C^W}$ is a partition of $W$ extending the corresponding partitions of $X$ and $Y$.
	Let $\cR$ denote the signature of $\Ef$. We shall use the notation introduced in the definition of diversification (in particular $\tilde R$ is the relation coming from $R$ by adding one more coordinate).
	We have to consider the following three cases.
	
	\paragraph{Case 1:} Both $X$ and $Y$ add $G$-orbits in $C^W$.
	
	In this case there is nothing to do, simply $W$ already carries a $\dup \Ef$-structure and the $G$-action satisfies (\ref{EqPiczfork}).
	
	\paragraph{Case 2:} One of the extensions adds a $G$-orbit in $C^W$ while the other one adds a $G$-orbit in $P^W$.
	
	We may assume $X = Z \cup A$, $Y = Z \cup B$ with $A \subs P^W \setminus Z$, $B \subs C^W \setminus Z$.
	Fix $b_0 \in B$.
	Then $P^Z$ has an $\Ef$-structure induced by $b_0$ in $Y$, because $P^Y = P^Z$.
	By the remarks above, as $\Ef \ne \sn \emptyset$, each model in $\Ef$ can be extended to an arbitrarily large finite model that is still in $\Ef$.
	So, let us choose some $\Ef$-structure on $P^W = P^Z \cup A$ extending the structure of $P^Z$ induced by $b_0$.
	Given $g \in G$, define the $\Ef$-structure on $P^W$ by using the $G$-action, namely,
	$$\tilde R(x_1, \dots, x_n, b_0^g) \iff \tilde R(x_1^{g^{-1}}, \dots, x_n^{g^{-1}}, b_0).$$
	By this way, the $G$-action becomes compatible with the $\Ef$-structure, namely, given $h \in G$ it holds that
	$$\tilde R(x_1^h, \dots, x_n^h, b_0^{gh}) \iff \tilde R(x_1^{g^{-1}}, \dots, x_n^{g^{-1}}, b_0) \iff \tilde R(x_1, \dots, x_n, b_0^g).$$
	
	\paragraph{Case 3:} Both $X$ and $Y$ add $G$-orbits in $P^W$.
	
	Here we essentially use the fact that $\Ef$ has strong amalgamations.
	Let $S \subs C^W = C^Z$ be a selector from the family of all $G$-orbits in $C^W$. Namely, $|S \cap O| = 1$ for every $G$-orbit $O \subs C^W$.
	Fix $s \in S$ and choose an $\Ef$-structure on $P^W = X^W \cup Y^W$ using the strong amalgamation property. This will be the structure induced by $s$.
	Next, given $g \in G$, define the $\Ef$-structure on $P^W$ induced by $s^g$ using the $G$-action, namely, for each relation $R$ define
	$$\tilde R(x_1, \dots, x_n, s^g) \iff \tilde R(x_1^{g^{-1}}, \dots, x_n^{g^{-1}}, s).$$
	In order to see that $W$ with such defined $\Ef$-structures is in $\dupg G \Ef$, it remains to check that condition (\ref{EqPiczfork}) holds.
	
	So fix an $n$-ary relation $R$ in the signature of $\Ef$ and fix $x_1, \dots, x_n \in P^W$, $y \in C^W$ such that $\tilde R(x_1, \dots, x_n, y)$ holds in $W$.
	As the $G$-action is free, there are uniquely determined $s \in S$ and $h \in G$ such that $y = s^h$.
	Fix $g \in G$.
	We have

	\begin{align*}
	\tilde R(x_1^g, \dots, x_n^g, y^g) &\iff \tilde R(x_1^g, \dots, x_n^g, s^{hg}) \iff \tilde R(x_1^{h^{-1}}, \dots, x_n^{h^{-1}}, s) \\
	&\iff \tilde R(x_1, \dots, x_n, s^h) \iff \tilde R(x_1,\dots, x_n, y).
	\end{align*}
	We conclude that $W \in \dupg G \Ef$, which completes the proof.
\end{pf}

\begin{lm}\label{LMerihgbwirug}
	Let $\Ef$ be a relational \fra\ class, let $X \in \dup \Ef$ be nonempty, and let $G$ be a finite group.
	Then there exists $X^G \in \dupg G \Ef$ such that $X$ embeds into $X^G$.
	Furthermore, if some $X_0 \subs X$ has a fixed free $G$-action, then at least one embedding of $X$ into $X^G$ preserves this action.
\end{lm}

Here, we consider $X^G$ as an $\Ef$-structure, forgetting the $G$-action.

\begin{pf}
	Let $X^G = \setof{x^g}{x \in X,\; g \in G}$, where we agree that $x^1 = x$ and $x^g \ne y^h$, unless $x = y$ and $g = h$.
	For each $c \in C^X$ extend the $\Ef$-structure induced by $c$ to 
	$(P^X)^G = \setof{x^g}{x \in P^X, \; g \in G}$.
	Finally, $c \in C^X$ and $g \in G$, define the $\Ef$-structure induced by $c^g$ by using the canonical $G$-action on $(P^X)^G$:
	$$\tilde R(x_1, \dots, x_n, y^g) \iff \tilde R(x_1^{g^{-1}}, \dots x_n^{g^{-1}}, y)$$
	for every $n$-ary relation $R$ in the signature of $\Ef$.
	As before, it is easy to check that the $G$-action preserves the $\dup \Ef$-structure, therefore $X^G \in \dupg G \Ef$.
	
	Finally, in case some $X_0$ has a fixed free $G$-action, we apply the procedure above to the set $X \setminus X_0$.
\end{pf}

\begin{thm}\label{ThmDiversitiesFree}
	Let $\Ef$ be a relational \fra\ class with strong amalgamations.
	Then $\dup \Ef$ is a \fra\ class and every finite group acts freely on its \fra\ limit.
\end{thm}

\begin{pf}
	By Lemma~\ref{LMKrucjalni}, $\dupg G \Ef$ and, in particular, $\dup \Ef = \dupg {\sn 1} \Ef$ is a \fra\ class.
	Let $U$ denote the \fra\ limit of $\Ef$ and let $\dupg G U$ denote the \fra\ limit of $\dupg G \Ef$.
	We write $\dup U$ for $\dupg {\sn1} U$. Since $\dupg G U$ has a canonical free $G$-action, we need to show that $\dup U$ is isomorphic to (the reduct of) $\dupg G U$ for every finite group $G$.
	For this aim, it suffices to check that $\dupg G U$ has the extension property with respect to $\dup \Ef$ (forgetting the $G$-action).
	This actually follows from Lemma~\ref{LMerihgbwirug}: Given a finite substructure $X_0 \subs \dupg G U$, given an embedding $\map{e}{X_0}{X}$ with $X \in \dup \Ef$, Lemma~\ref{LMerihgbwirug} provides a $\dupg G \Ef$-structure on $X$ extending the $\dupg G \Ef$-structure of $X_0$, therefore there is an embedding $\map f X {\dupg G U}$ such that $f e$ is the inclusion $X_0 \subs U$.
	Forgetting the group action, we see that $\dupg G U$ has the extension property with respect to $\dup \Ef$, therefore it is the \fra\ limit of $\dup \Ef$.
\end{pf}

\subsection{An application: consumer-product models}

Let $\procons$ be the class of all finite models $\seq{M, P^M, C^M, L^M}$, where $P, C$ are unary predicates, $L$ is a ternary relation, and the following two axioms are satisfied.
\begin{enumerate}[itemsep=0pt]
	\item[(CP1)] $\dn{P^M}{C^M}$ is a partition of $M$.
	\item[(CP2)] For each $c \in C^M$, the relation $L^M(\argum, \argum, c)$ is a strict linear ordering of $P^M$.
\end{enumerate}
We shall rather write $x <_c y$ instead of $L^M(x,y,c)$.
The elements of $P^M$ will be called \emph{products} while the elements of $C^M$ will be called \emph{consumers}.
The idea is that $P^M$ consists of certain items (goods) of the same type. Consumers have their personal preferences so that, given two different products, each individual consumer can say which one is better or more desirable for her/him, thus implicitly defining a linear ordering on the set of items.
Clearly, $\procons = \dup \El$, the diversification of the \fra\ class of all finite linearly ordered sets.
Structures in the class $\loc \procons$ will be called \emph{consumer-product models}.
The infinite ones may represent limits of some evolution processes where either the number of products or the number of consumers or both tend to infinity.

By Theorem~\ref{ThmDiversitiesFree}, $\procons$ is a \fra\ class. Let $\U$ denote its \fra\ limit.
Also by Theorem~\ref{ThmDiversitiesFree}, $\aut(\U)$ contains copies of all finite groups, which proves part of Theorem~\ref{ThmMejnOne}.

It is rather obvious that $\U$ has infinitely many consumers, infinitely many products, and for each consumer $c$ the relation $<_c$ defines an ordering of $P^\U$ isomorphic to $\pair{\Qyu}{<}$.
Furthermore, no two consumers have the same preferences, namely, if $c \ne d$ are in $C^\U$ then there are products $p, q$ such that $p <_c q$ and $q <_d p$.
This is a straightforward consequence of the extension property. Namely, $N = \dn c d$ with $P^N=\emptyset$ is a submodel of $\U$ and $M \sups N$ defined by $P^M = \dn p q$ and $p <_c q$, $q <_d p$ is a consumer-product model that has to be realized inside $\U$.

Fix a consumer $c \in \U$ and define
$$K_c = \setof{h \in \aut(\U)}{h(c) = c}.$$
If $h \in K_c$ and $h \rest P^\U$ is identity, then $h = \id \U$, because of the remark above.
Thus, the mapping $h \mapsto h \rest P^\U$, defined on $K_c$ into the group of permutations of $P^\U$, has trivial kernel. This shows that $K_c$ is isomorphic to a subgroup of $\aut(\Qyu,<)$, therefore it is torsion-free.
On the other hand, $K_c$ has a countable index in $\aut(\U)$, because $C^\U$ is countable and each coset of $K_c$ is defined by some (any) automorphism moving $c$ to another consumer $d \in C^\U$.
More precisely, if $f_0, f_1 \in \aut(\U)$ are such that $f_i(c) = d$ for $i=0,1$, then $f_1^{-1} \cmp f_0 \in K_c$.
We conclude that $S_\infty$ does not embed into $\aut(\U)$, by Corollary~\ref{CorWNebwgibEIG}.
On the other hand, $\loc \procons$ has infinite models without products or without consumers, whose automorphism groups are clearly isomorphic to $S_\infty$.
This shows that $\aut(\U)$ is not universal, which completes the proof of Theorem~\ref{ThmMejnOne}.

Let us finally note that while $\U$ has an involution, since $\Zee_2$ acts freely on it, there exist involutions of its finite submodels such that all their extensions have infinite order.
Namely, fix a finite $S \subs \U$ consisting of consumers only. Then $\aut(S)$ is the group of all permutations of $S$. Assume $|S|>2$ and let $h \in \aut(S)$ be a nontrivial involution with a fixed point $c \in S$. By the arguments above, if $\tilde h \in \aut(\U)$ extends $h$ then its restriction on the set of products is a non-trivial automorphism of $\pair{P^\U}{<_c} \iso \pair \Qyu <$, therefore it has an infinite order.
This is a clear evidence that $\procons$ cannot have a Kat\v etov functor (cf.~\cite{KubMas}) and its \fra\ limit is not uniformly homogeneous.

\section{Mixed sums}

We now describe another operation, this time on pairs of \fra\ classes, involving bipartite graph structures.

Let $\Ef$, $\Gee$ be two classes of models (possibly of different signatures). We define their \emph{mixed sum} $\Ef \misa \Gee$ to be the class of all structures of the form $M = L^M \cup R^M$, where $L^M \cap R^M = \emptyset$, $L^M \in \Ef$, $R^M \in \Gee$, and additionally $\triple{L^M}{R^M}{\sim^M}$ is a bipartite graph. More precisely, $\pair M{\sim^M}$ is a graph such that $x \sim^M y$ implies $x \in L^M$, $y \in R^M$ or vice versa.
Formally we should assume that the signatures of $\Ef$ and $\Gee$ are disjoint\footnote{As this concerns symbols only, there is no problem here.} and the signature of $\Ef \misa \Gee$ is their union, together with three new predicates: $L, R, \sim$.

\begin{lm}
	Assume $\Ef$, $\Gee$ are relational \fra\ classes, each of them having strong amalgamations.
	Then $\Ef \misa \Gee$ is a \fra\ class with strong amalgamations.
\end{lm}

\begin{pf}
	It is clear that $\Ef \misa \Gee$ is hereditary and has countably many types.
	Let $\map f Z X$, $\map g Z Y$ be two embeddings with $Z, X, Y \in \Ef \misa \Gee$.
	We may assume that $f$, $g$ are inclusions and $Z = X \cap Y$. It suffices to check that $W := X \cup Y$ carries a structure of $\Ef \misa \Gee$.
	
	Using the strong amalgamation property of $\Ef$, we find $L^W \in \Ef$ so that $L^X \cup L^Y \subs L^W$ and the inclusions are embeddings. In fact, we may assume (although it is irrelevant here) that $L^W = L^X \cup L^Y$, because $\Ef$ is relational.
	We do the same with $R^X$, $R^Y$, obtaining $R^W \in \Gee$ containing $R^X \cup R^Y$.
	Let $W = L^W \cup R^W$ (of course, we assume that $L^W \cap R^W = \emptyset$).
	Finally, we define $x \sim^W y$ if and only if $x \sim^X y$ or $x \sim^Y y$. Now $W \in \Ef \misa \Gee$ and the inclusions $X \subs W$, $Y \subs W$ are embeddings.
	This shows that $\Ef \misa \Gee$ has strong amalgamations. The joint embedding property follows from the amalgamation property, because the empty set is a model both in $\Ef$ and in $\Gee$.
\end{pf}

Note that once we allow functions in the signatures,  $\Ef \misa \Gee$ still has the amalgamation property, as long as $\Ef$, $\Gee$ have strong amalgamations and all the models are finite (if some $X \in \Ef \misa \Gee$ is infinite then there are uncountably many bipartite graph structures on $\pair{L^X}{R^X}$). On the other hand, if there are some constants in the languages of $\Ef$ and $\Gee$ then the class $\Ef \misa \Gee$ fails the joint embedding property.

We say that a class of models $\Ef$ is \emph{non-degenrate} if $\Ef \ne \sn \emptyset$.

\begin{tw}\label{THMsbgisbgse}
	Assume $\Ef$, $\Gee$ are non-degenerate relational \fra\ classes, both with the strong amalgamation property. Let $\U_\Ef$, $\U_\Gee$ denote their \fra\ limits, and let $\U$ denote the \fra\ limit of $\Ef \misa \Gee$.
	
	Then $L^\U \iso \U_\Ef$, $R^\U \iso \U_\Gee$ and $\triple{L^\U}{R^\U}{\sim^\U}$ is the universal homogeneous bipartite graph.
	Furthermore, the restriction mappings $h \mapsto h \rest L^\U$ and $h \mapsto h \rest R^\U$ are embeddings of $\aut(\U)$ into $\aut(\U_\Ef)$ and $\aut(\U_\Gee)$, respectively.
\end{tw}

\begin{pf}
	First of all, notice that  $\U_\Ef$ and $\U_\Gee$ are infinite, because of the strong amalgamation property and the existence of non-empty models.
	
	In order to show that $L^\U \iso \U_\Ef$, we check the extension property.
	Fix $A \subs B \in \Ef$ and an embedding $\map e A {L^\U}$, where $L^\U$ is viewed as a structure in $\Ef$. Modifying the language, $B$ can be regarded as a structure in $\Ef \misa \Gee$, where $L^B = B$, $R^B = \emptyset$, and the bipartite graph relation is empty.
	Now, using the extension property of $\U$, there is an embedding $\map f B \U$ extending $e$.
	Coming back to the original language of $B$, we conclude that $f$ is an embedding of structures in $\Ef$. Hence $L^\U \iso \U_\Ef$.
	The same arguments show that $R^\U \iso \U_\Gee$.

	We now check, using Proposition~\ref{PROPbhvbivbw}, that $\triple{L^\U}{R^\U}{\sim^\U}$ is the universal homogeneous bipartite graph.
	Again, fix a bipartite graph $B = \triple{L^B}{R^B}{\sim^B}$, its subgraph $A$, and an embedding $\map e A \U$, where now $\U$ is viewed as the bipartite graph (so the signature consists of three symbols: $L$, $R$, $\sim$).
	Using $e^{-1}$, we endow $A$ with the $\Ef \misa \Gee$-structure, so that now $A \in \Ef \misa \Gee$.
	Using the fact that both $\Ef$ and $\Gee$ have arbitrarily large models, we may find an $\Ef \misa \Gee$-structure on $B$ extending that of $A$.
	Now $e$ is an embedding of $\Ef \misa \Gee$-structures and, by the extension property of $\U$, we find an embedding $\map f B \U$ satisfying $f \rest A = B$.
	Forgetting the $\Ef \misa \Gee$-structure, leaving only the bipartite graph relation, we conclude that $f$ is an embedding of bipartite graphs, showing that $\triple \U {L^\U} {R^\U}$ is the \fra\ limit of finite bipartite graphs.

	Finally, fix $h \in \aut(\U)$ and suppose $h \ne \id \U$.
	Then $h(b_0) \ne b_0$ for some $b_0 \in R^\U$.
	Let $A \subs \U$ be any finite model containg $\dn{b_0}{h(b_0)}$. Extend $A$ to a model $B \in \Ef \misa \Gee$ so that for some $a_0 \in L^B \setminus L^A$ we have that $a_0 \sim^B b_0$ and $a_0 \not \sim h(b_0)$. This is possible, because $\Ef$ has arbitrarily large models and there are no restrictions for the bipartite graph relation.
	By the extension property, $B$ is realized in $\U$, so we may assume $a_0 \in L^\U$. It cannot be the case that $h(a_0) = a_0$, therefore $h \rest L^\U$ is not identity.
	This shows that the restriction map $h \mapsto h \rest L^\U$ has a trivial kernel. The mixed sum is symmetric, therefore the same arguments apply to the map $h \mapsto h \rest R^\U$.
\end{pf}

As an application, consider $\Ef$ to be the class of all finite sets (so the signature of $\Ef$ is empty) and $\Gee$ to be the class of all finite linearly ordered sets. Let $G = \aut(\U)$, where $\U$ is the \fra\ limit of $\Ef \misa \Gee$.
By Theorem~\ref{THMsbgisbgse}, $G$ embeds into $\aut(\Qyu, <)$, therefore it is torsion-free.
On the other hand, among models of $\Ef \misa \Gee$ we may find those having the empty bipartite graph relation, therefore all finite symmetric groups appear as $\aut(X)$ with $X \in \Ef \misa \Gee$.
This proves Theorem~\ref{ThmMejnTwo}.

\section{Rotating machines}

We now present another example of a \fra\ class of finite models with non-trivial abelian automorphism groups, whereas the automorphism group of its \fra\ limit is torsion-free. This will prove Theorem~\ref{ThmMejnTwo}.

Consider the signature $\El$ with two binary predicates $<$, $\sim$ and one unary function symbol $\suk{}$. In a model of this language, the mapping $x \mapsto \suk x$ will be called the \emph{successor operation} whenever it is one-to-one.
A \emph{rotating wheel} is a finite structure $M$ of this language, satisfying the following conditions:
\begin{enumerate}[itemsep=0pt]
	\item[(W0)] The successor operation $x \mapsto \suk x$ is bijective and has exactly one orbit (in other words, the whole set $M$ forms a cycle with respect to this operation).
	\item[(W1)] The relations $<$, $\sim$ are empty.
\end{enumerate}
Note that every automorphism of a rotating wheel $M$ is actually a power of the successor operation.
In other words, $\aut(M)$ is cyclic of order $|M|$, generated by the successor operation.
Note also that every rotating wheel is isomorphic to $\pair{\Zee_n}{\suk{}}$, where $\suk{x} = x +_n 1$, where $+_n$ denotes the addition modulo $n$.

A \emph{rotating machine} is a finite model $M$ of the signature $\El$ satisfying the following axioms.
\begin{enumerate}[itemsep=0pt]
	\item[(V0)] $M$ is a disjoint union of rotating wheels.
	\item[(V1)] $<$ is a strict partial order and $\sim$ is an undirected graph relation.
	\item[(V2)] If $x < y$ or $x \sim y$ then $x$ and $y$ belong to different rotating wheels.
	\item[(V3)]	If $C$ and $D$ are different rotating wheels then either $C < D$ or $D < C$.
	\item[(V4)] $\sim$ is compatible with the successor operation, namely, $x \sim y \implies \suk x \sim \suk y$.
\end{enumerate}
Concerning (V3), we use the abbreviation $A < B$ meaning $a < b$ for every $a \in A$, $b \in B$.
Condition (V4) is crucial, it says that different rotating wheels are connected in the sense that rotating one of them induces a suitable rotation of the other. We will use it in the proof of Lemma~\ref{LMsdbguohwe} below.

Let $\Vee$ denote the class of all rotating machines.

\begin{lm}
	$\Vee$ is a \fra\ class of finite models.
\end{lm}

\begin{pf}
	Rather trivial, the amalgamation property can be proved almost in the same way as for the linearly ordered graphs.
\end{pf}

Let $\bW$ denote the \fra\ limit of $\Vee$.
\begin{lm}\label{LMsdbguohwe}
	The group $\aut(\bW)$ is torsion-free.
\end{lm}

\begin{pf}
	Fix $h \in \aut(\Wu)$, $h \ne \id{\Wu}$ and suppose $h$ has a finite order $k > 1$.
	Note that each rotating wheel must be invariant under $h$, since otherwise $h$ would induce a non-trivial isomorphism of the linearly ordered set of all rotating wheels, therefore $h$ would have infinite order.
	
	There exists a rotating wheel $C \subs \Wu$ such that $h \rest C$ has order $k$. Indeed, if $\ell$ is the maximum of the orders of $h$ restricted to rotating wheels, then $\ell \loe k$ and $h^\ell = \id{\Wu}$, therefore $\ell = k$.
	
	Let $n = |C|$, so $C$ is isomorphic to $\Zee_n$ with the standard successor operation. Let us use the enumeration $C = \{0^C, \dots, (n-1)^C\}$, where $\suk{i^C}=(i +_n 1)^C$ for $i < n$.
	Let $D = \Zee_m$ with $m = k n$ and define the graph relation between $C$ and $D$ by
	$$x^C \sim y^D \iff x = y \mod n.$$
	We use the enumeration $D = \{0^D, \dots, (m-1)^D\}$.
	Of course, $\sim$ must be symmetric, therefore we also define $y^D \sim x^C \iff x^C \sim y^D$. We define $<$ so that $C < D$ (actually the ordering plays no role here).
	We need to check that $C \cup D$ is a rotating machine and the only possible obstacle is condition (V4).
	
	Fix $x, y$ such that $x^C \sim y^D$.
	If $x < n-1$ and $y < m-1$ then clearly $\suk {x^C} \sim \suk {y^D}$.
	Suppose $x = n-1$. Then $\suk{x^C} = 0^C$ and $y = n y' - 1$ for some integer $y'$; hence $y +_m 1$ is divisible by $n$, as it is either $y + 1$ or $0$. Thus in this case $\suk{(x^C)} \sim \suk{y^D}$.
	Finally, suppose $y = m-1$.
	Knowing that $n$ divides $m$, we see that necessarily $x = n-1$, therefore by the previous case $\suk{x^C} \sim \suk{y^D}$.
	
	We have verified condition (V4), concluding that $C \cup D$ is indeed a rotating machine.
	Thus, we may assume that it is already contained in $\Wu$.
	
	Let $h(0^C) = a^C$, $a \in \Zee_n$ and let $h(0^D) = b^D$, $b \in \Zee_m$.
	Then $a^C \sim b^D$, because $0^C \sim 0^D$.
	Note that $h^k(0^C) = k a \mod n$ and $h^k(0^D) = k b \mod m$.
	Thus, there is $c \in \nat$ with $k b = c m = c k n$. Hence $b = c n$. It follows that $a$ is divisible by $n$, which is possible only if $a=0$, a contradiction, because $h \rest C$ is not the identity.
\end{pf}

Finally, among rotating machines we may find infinitely many with empty edge relation and their automorphism groups are arbitrary large finite products of cyclic groups.
Furthermore, given a rotating machine, its automorphism group embeds into the product of all automorphism groups of its rotating wheels, which in turn are finite cyclic groups. In particular, for every rotating machine $M$, the group $\aut(M)$ is abelian.
This proves Theorem~\ref{ThmMejnTwoAnAhalf}.

\paragraph{Acknowledgments.}
The authors would like to thank Alex Kruckman for useful comments and for pointing out reference~\cite{Baudisch}.

\end{document}